\documentclass{article}
\usepackage[T2A]{fontenc}
\usepackage[cp866]{inputenc}
\usepackage[english]{babel}
\usepackage[tbtags]{amsmath}
\usepackage{amsfonts,amssymb,mathrsfs,amscd,longtable, accents}
\usepackage{mathrsfs}

\usepackage{hypbmsec}
\usepackage{mathtext}

\newtheorem{theorem}{Theorem}

\let\epsilon\varepsilon

\begin{document}
\title{Linear Lie Groups All of Whose Irreducible\\Finite-Dimensional Not Necessarily
Unitary\\Representations Are of Bounded Dimension\\and Separate the Points of
the Group}

\author{A.~I.~Shtern}

\begin{date}{10.10.2025}\end{date}

 \maketitle

 \centerline{Lomonosov Moscow State University,}
 \centerline{Department of Mechanics and Mathematics;}
 \smallskip
 \centerline{NRC ``Kurchatov Institute'' - SRISA}
 \smallskip
 \centerline{E-mail: \tt{rroww\@mail.ru}, \tt{aishtern\@mtu-net.ru}}
 \smallskip
 \centerline{2020 \textit{Mathematics Subject Classification.}
 Primary 22A99, Secondary 22A25 22E25.}

\begin{abstract}
 We prove that all linear Lie groups satisfying the conditions listed in the
title are finite extensions of commutative Lie groups.
\end{abstract}

 \centerline{\textbf{Keywords.} Lie groups, finite-dimensionally representable
 topological group.}

\markboth{Shtern}{Linear Lie Groups with Finite-Dimensional Representations}

\section{Introduction}
\label{s1}

In this note, we prove that every Lie group satisfying the conditions listed in
the title is a finite extension of a commutative Lie group admitting the set of
continuous characters that separates the elements of the group.

\section{Preliminaries}\label{s2}

Recall that a locally compact group (in particular, a Lie group) all of whose
continuous unitary representations are finite-dimensional and their dimensions
are jointly bounded is a finite extension of a commutative locally compact
(Lie) group~\cite{1},~\cite{2}.

 \section{Main Results}\label{s3}

\begin{theorem}\label{t1}
  A linear Lie group for which all its irreducible continuous finite-dimensional not
necessarily unitary linear representations are of bounded degree and separate
the points of the group is a finite extension of a commutative Lie group.
\end{theorem}

\textbf{Proof.} \textrm{Let $G$ be a Lie group satisfying the conditions of the
theorem.}

\textrm{Let $G_0$ be the connected component of the identity of~$G$. The group $G_0$
obviously satisfies the conditions of the theorem. Let $$G_0=LR$$ be a Levi
decomposition of~$G$, where $L$ is a semisimple subgroup of~$G_0$ and~$R$ is
the radical of~$G_0$.}

\textrm{Then $L$ is the identity group, because otherwise its adjoint group is a
quotient group of~$G_0$ which is a nontrivial semisimple Lie group, and hence
has continuous finite-dimensional representations of arbitrary large
finite dimensions, which contradicts the assumption.}

\textrm{Therefore, $G_0$ is solvable. A continuous irreducible representation of a 
connected solvable Lie group is one-dimensional by the Lie theorem. Hence,
all irreducible representations of~$G_0$ are one-dimensional. If the 
connected solvable group $G_0$ is not commutative, then $G_0$ has finite-dimensional
representations, which contradicts the conditions of the theorem that hold
for~$G_0$. Thus, $G_0$ is commutative.}

\textrm{The quotient group $G/G_0$ is discrete and satisfies the conditions of the
theorem. Then the group algebra $l^1(G/G_0)$ satisfies the identity (1)
of~Subsec.~3.6.1 in~\cite{3} for $r=r(n)$, where the dimensions of the
representations mentioned in the condition of the theorem do not exceed $n$ and
$r(n)$ is the least value of $r$ for which this identity holds.}

\textrm{By the Gel'fand--Raikov theorem, this group has a set of irreducible
representations in Hilbert spaces that separates the points of~$G/G_0$.}

\textrm{As was in fact proved in~Proposition~3.6.3 of~\cite{3}, the dimensions of these
representations are bounded by~$n$ and, by Thoma's theorem~\cite{4}, $G/G_0$ is
a finite extension of a commutative normal subgroup $M$ of~$G/G_0$.}

\textrm{The group~$M$ is discrete, since~$G/G_0$ is, and hence the (commutative)
closure of~$M$ with respect to the topology of~$G$ coincides with~$M$.}

\textrm{Since the extension $H$ of $G_0$ by $M$ is a (disconnected) Lie group, it
follows that a theorem of Calvin Moore~\cite{5} (concerning locally compact
groups all of whose continuous unitary representations are finite-dimensional
and have bounded dimensions) holds for~$H$.}

\textrm{Hence, if the extension $H$ of $G_0$ by $M$ is not a finite extension of a
commutative Lie group~$G_0$, then the family of irreducible finite-dimensional continuous
unitary representations of $H$ either does not have bounded dimensions or does not 
separate the points of the group, which contradicts the conditions of the theorem for~$G$, 
since they obviously remain valid for the subgroup~$H$ of~$G$.}

\textrm{Thus, $H$ is a finite extension of a commutative Lie group and, therefore, so
is~$G$.}

\textrm{Since the maximal commutative Lie subgroup of~$H$ also has the properties
imposed on~$G$, it follows that this commutative group is finite-dimensionally
representable, and hence has the set of continuous characters separating the
elements of the group.}

\textrm{This completes the proof of the theorem.}
\smallskip

\section{Discussion}\label{s4}

In~\cite{5}, we proved a theorem for Hausdorff topological groups all of whose
irreducible unitary finite-dimensional representations are finite-dimensional
and their dimensions are jointly bounded for which the family of these
finite-dimensional representations separates the points of the group: the
theorem claims that every group of this kind is a finite extension of a
commutative topological group whose continuous characters separate the points
of the group.

\section*{Funding}

The research was supported by the Scientific Research Institute for System
Analysis of the National Research Centre ``Kurchatov Institute'' according to
the project FNEF-2024-0001.

\end{document}